\documentclass[10pt, leqno]{amsart}
\usepackage{amscd}
\usepackage{amssymb}
\usepackage{amsmath}
\usepackage{amsthm}
\usepackage[cmtip,all]{xypic}
\usepackage{pifont}

\theoremstyle{plain}
\newtheorem{Lem}{Lemma}[section]

\newtheorem{Thm}[Lem]{Theorem}

{\theoremstyle{definition} 

\newtheorem{Rk}[Lem]{Remark}
\newtheorem{Def}[Lem]{Definition}}

\newcommand{\lhat}{\hat{L}}
\newcommand{\zig}{\addtocounter{Lem}{1}\tag{\theLem}}

\def\:{\colon}

\begin{document}
\title{Rognes's theory of Galois extensions and the continuous action of $G_n$ on $E_n$}
\author{Daniel G. Davis}
\date{May 14, 2004; abstract written on 11/30/06}
\begin{abstract}
Let us take for granted that $L_{K(n)}S^0 \rightarrow E_n$ is 
some kind of a $G_n$-Galois extension. Of course, this is in the setting of 
continuous $G_n$-spectra. How much structure does this continuous $G$-Galois 
extension have? How much structure does one want to build into this 
notion to obtain useful conclusions? If the author's conjecture that ``$E_n/I$, for a cofinal collection of $I$'s, is a discrete $G_n$-symmetric ring 
spectrum" is true, what additional structure does this give the 
continuous $G_n$-Galois extension? Is it useful or merely beautiful? This 
paper is an exploration of how to answer these questions. This preprint 
arose as a letter to John Rognes, whom he thanks for a helpful 
conversation in Rosendal. This paper was written before John's preprints 
(the initial version and the final one) on Galois extensions were 
available. 
\end{abstract}
\maketitle
\begin{section}{Introduction}
In recent years, John Rognes has given various talks introducing his theory of Galois extensions for commutative $S$-algebras, and several manuscripts about this topic are available at his website. One family of Galois extensions that he has discussed are those that arise from the action of the extended Morava stabilizer group $G_n=S_n \rtimes \mathrm{Gal}(\mathbb{F}_{p^n}/{\mathbb{F}_p})$ on the Lubin-Tate spectrum $E_n$ through automorphisms of commutative $S$-algebras. In considering this new theory, the author has found these $E_n$-related Galois extensions quite interesting. 
\par
In the author's 2003 thesis (\cite{thesis} - see \cite{note} for a short summary), for a profinite group $G$, he developed (making explicit ideas that were implicit in the literature, especially in the foundational work of Bob Thomason and Rick Jardine, and work by Paul Goerss and Steve Mitchell) the notions of continuous action and homotopy fixed points for discrete $G$-spectra and their towers. For a closed subgroup $G$ in $G_n$, Ethan Devinatz and Mike Hopkins defined $$E_n^{hG} := L_{K(n)}(\mathrm{colim} \, _i \, E_n^{hU_iG})$$ (see \cite[Def. 1.5]{DH}). The author applied their work \cite{DH} and this definition to show that $E_n$ is a continuous $G_n$-spectrum with homotopy fixed point spectra, defined using the continuous action, for closed subgroups $G$ in $G_n$. 
\par
The above formula for $E_n^{hG}$ follows a convention that is used throughout this paper: $E_n^{hG}$ is used to denote  both the homotopy fixed point spectra of Devinatz and Hopkins, and the homotopy fixed point spectra defined with respect to the continuous action of $G$ on $E_n$ (defined by the author in \cite{thesis}), since the author showed that these constructions are isomorphic in the stable homotopy category (in \cite{thesis}). 
\par
Since the author does not possess a detailed account of Rognes's ideas (nor has he had the fortune of hearing Rognes speak), and, believing that the machinery of his thesis could be useful for the theory of Galois extensions, the author wrote this paper, to help him precisely understand Rognes's examples and to see more clearly exactly what kind of Galois extensions arise in Lubin-Tate theory. Thus, the job of this paper is primarily to study the extensions that arise from $E_n$, with its $G_n$-action, and to consider what kinds of definitions of Galois extension are needed to fit the actual structures. 
\par
Part of our work in this paper depends on results that are not yet known to be true. Beginning in \S 3, we assume that the discrete $G_n$-module $\pi_\ast(E_n \wedge M_I)$ can be realized by a spectrum, abusively labeled $F_n \wedge M_I$, that is a discrete $G_n$-symmetric ring spectrum, that is, a discrete $G_n$-spectrum that is a ring object in the category of symmetric spectra, and whose discrete $G_n$-action is by ring maps. We also assume the existence of certain model categories for various categories of discrete $G$-symmetric spectra. See Remark \ref{unproven} for more details about our assumptions. The author hopes and believes that these assumptions are true.
\par
If we only assume what is already known and if we modify Definition \ref{pro} (for ``$K(n)$-local $G$-pro-Galois extension'') in an obvious way, then our main results, stated as (4) and (6) below, are still true. We make the above assumptions because they allow us to paint a more coherent and structured picture than would be possible otherwise, and because we hope that the picture will turn out to be correct. 
\vspace{.05in}
\par \noindent
\textbf{Summary of Main Results.}
To ease the notation, we write $\lhat$ in place of $L_{K(n)}$. We summarize the types of extensions and examples that are considered in this paper by listing our main results (given the above assumptions):
\begin{enumerate}
\item
Theorem \ref{ex1}: Given an open normal subgroup $U$ of $G_n$, the map of commutative $S$-algebras $\alpha(U) \: \lhat S^0 \rightarrow E_n^{hU}$ is a $K(n)$-local $G_n/U$-Galois extension.
\item
Theorem \ref{assoc}: The map $\gamma(U,I) \: \lhat M_I \rightarrow (F_n \wedge M_I)^{hU} \cong E_n^{hU} \wedge M_I$ is an associative $G_n/U$-Galois extension.
\item
Theorem \ref{thm}: The map $\gamma(I) = \mathrm{colim} \, _i \, \gamma(U_i,I) \: \lhat M_I \rightarrow F_n \wedge M_I$ is an associative filtered $G_n$-Galois extension.
\item
Theorem \ref{ex4}: The map $\mathrm{holim} \, _I \, \gamma(I) \: \lhat S^0 \rightarrow E_n$ is a $K(n)$-local $G_n$-pro-Galois extension. Also, we explain why $\lhat S^0 \rightarrow E_n$ is a strongly $K(n)$-local filtered $G_n$-pro-Galois extension.  
\item
Theorem \ref{ex5}: Given any closed subgroup $G$ of $G_n$, the directed system $\{\beta(G,i,I) \: E_n^{hU_iG} \wedge M_I \rightarrow E_n^{hU_i} \wedge M_I\}_i$, of associative $U_iG/{U_i}$-Galois extensions, makes the map $$\ \ \ \ \ \beta(G,I) = \mathrm{colim} \, _i \, \beta(G,i,I) \: (E_n^{hG} \wedge M_I) \rightarrow F_n \wedge M_I \simeq E_n \wedge M_I$$ an associative filtered $G$-Galois extension.
\item
Theorem \ref{ex6}: The inverse system $\{\beta(G,I)\}_I$ of $K(n)$-local associative $G$-Galois extensions makes
the map $\beta(G) = \mathrm{holim} \, _I \, \beta(G,I) \: E_n^{hG} \rightarrow E_n$ a $K(n)$-local $G$-pro-Galois extension.
\item
Theorem \ref{ex7}: The map $\alpha(U) = \mathrm{holim} \, _I \, \gamma(U,I)$ is a $K(n)$-local $G_n/U$-pro-Galois extension.
\end{enumerate}
\vspace{.05in}
\par
\noindent
\textbf{Notation and Conventions.} Often, when we use results from \cite{DH} and \cite{thesis}, we do not give references. Throughout this paper, $U$ is an open normal subgroup of $G_n$. $Sp$ is the model category $\mathrm{(spectra)}^{\mathrm{stable}}$ of Bousfield-Friedlander spectra. We often use the symbol $\cong$ to denote isomorphism in the stable homotopy category. Whenever necessary, we assume that our commutative $S$-algebras are cell commutative, and, given an $S$-algebra $R$, that our $R$-modules are cell $R$-modules. All colimits of spectra are formed in $S$-modules, $Sp$, or $Sp^\Sigma$, the model category of symmetric spectra of simplicial sets; which category is used will be clear from context. Whenever necessary, we view an $S$-module as a (symmetric) spectrum of simplicial sets, and vice versa.
\par
Given a profinite group $G$, if a colimit or limit is indexed by a collection $\{N\}$, then $\{N\}$ is a cofinal collection of open normal subgroups of $G$. Also, $Sp_{_G}$ is the category of discrete $G$-spectra, and, given $X \in Sp_{_G}$, $X_{f,G}$ denotes the spectrum obtained from factoring $X \rightarrow \ast$ as $X \rightarrow X_{f,G} \rightarrow \ast$, a trivial cofibration, followed by a fibration, in $Sp_{_G}$. Then, by definition, $X^{hG} = (X_{f,G})^G.$
\par
If $X$ is a (pointed) discrete $G$-set, simplicial set, or spectrum, then $\mathrm{Map}_c(G,X)$ is the (pointed) discrete $G$-set, simplicial set, or spectrum of continuous maps from $G$ to $X$, and $\Gamma_G^\bullet (X)$ is the canonical cosimplicial (pointed) discrete $G$-set, simplicial set, or spectrum determined by the triple that is formed from $\mathrm{Map}_c(G,-)$. We will use the fact that if the profinite group $G$ has finite virtual cohomological dimension, then $\mathrm{holim} \, _\Delta \, \bigl(\Gamma^\bullet_G (X_{f,G})\bigr)^G$ is a model for $X^{hG}$ (see \cite{thesis}).
\vspace{.05in}
\par \noindent
\textbf{Note to Reader:} The author wants to point out that besides the assumption of the validity of certain conjectural remarks, discussed above, the careful reader will notice that there are several other technical problems with this paper, which we now list. (a) We go back and forth between $S$-modules and spectra of simplicial sets frequently, and this movement is less than desirable. An ideal presentation of the various Galois extensions considered here would do everything in a single category of spectra. (b) Since the colimit in $S$-modules of $S$-algebras is not, in general, the colimit in the category of $S$-algebras (see \cite[II, Prop. 7.4]{EKMM}), the colimit in Definition \ref{galfil} should be handled in a better way. (c) At places where we would like to have point-set level maps of (commutative) $S$-algebras that are weak equivalences, we often have only isomorphisms in the stable homotopy category. Though the author has not ironed out these technicalities, he still believes that the ideas in this paper are essentially correct and worthwhile.  
\vspace{.05in}  
\par
\noindent
\textbf{Relationship to Work of Rognes.} The inspiration for this paper is the work of John Rognes. Some of the definitions and examples are originally due to him, and those that are not were motivated by his work. At the time of writing, the author does not know how much of the details of what he has written here is already known to Rognes. To make the relationship between this work and Rognes's clearer, the author notes the following. (a) Definition \ref{gal} is from \cite[pg. 1]{Rognes2000}. (b) Definition \ref{gallocal} closely follows Rognes's definition of an $E$-local $G$-Galois extension (see \cite{Rognes2003}). (c) Definition \ref{galfil} was motivated by Rognes's notion of a $G$-pro-Galois extension (see \cite[pg. 1]{Rognes2000}). (d) Our main theorem, Theorem \ref{ex4}, is a slight modification of a result due to Rognes (see \cite[pg. 2, (6)]{Rognes2000}, \cite[pg. 6, (4)]{Rognes2002}, and \cite{Rognes2003}). The author readily acknowledges that he has perhaps omitted ways that this work is already known by, written up by, or indebted to Rognes. 
\vspace{.05in}
\par
\noindent
\textbf{Acknowledgements.} First of all, I thank Jeff Smith for his encouragement in my work on continuous $G$-spectra in chromatic stable homotopy theory. I thank John for sending me \cite{Rognes2002}, which was helpful. I thank Paul Goerss for many helpful conversations in the past several years about homotopy fixed points and various categories of spectra, and for discussions about the idea that $\pi_\ast(E_n \wedge M_I)$ can be realized by a discrete $G_n$-symmetric ring spectrum. Also, I thank Paul, Jim McClure, and Clarence Wilkerson for their encouragement.
\end{section}
\begin{section}{Galois extensions for finite groups}
\begin{Def}\label{gal}
Let $G$ be a finite group. A map $A \rightarrow B$ of commutative $S$-algebras is a {\em $G$-Galois extension} if the following conditions hold:
\begin{enumerate}
\item
$G$ acts on $B$ through $A$-algebra maps.
\item
The natural map $A \rightarrow B^{hG}$ is a weak equivalence. 
\item
There is a weak equivalence $B \wedge_A B \simeq \prod_G \, B,$ where $B$ is an $A$-module and $-\wedge_A -$ is the smash product in the category of $A$-modules.
\end{enumerate}
\end{Def}
\begin{Rk}
The theorem below is due to Rognes (\cite{Rognes2002}, \cite{Rognes2003}). Since the author knows of no written proof, he attempts one. The proof below is complete, except for the unjustified step marked with a question mark. The unjustified step is asking for a holim to commute with a colimit (specifically, a coequalizer). The author would like to know how to complete the proof, and then prove the same result when $G$ is profinite, in which case the result could be used with Galois extensions that are considered later in this paper.
\end{Rk}
\begin{Thm}\label{unsure}
Let $G$ be a finite group and let $N$ be a normal subgroup. If $A \rightarrow B$ is a $G$-Galois extension, then $A \rightarrow B^{hN}$ is a $G/N$-Galois extension.
\end{Thm}
\begin{proof}
We know that $B^{hG} \simeq A$ and $B \wedge_A B \simeq \prod_G \, B.$ Regarding $G$ as a profinite group, with $N$ an open subgroup, $(B^{hN})^{hG/N} \simeq B^{hG} \simeq A$ (see Lemma \ref{jardine}). It remains only to show that $B^{hN} \wedge_A B^{hN} \simeq \prod \, _{G/N} \, B^{hN}.$ We have (where $N$ acts on $B$ and $\mathbf{N}$ acts on $\mathbf{B}$):
\begin{align*}
B^{hN} \wedge_A B^{hN} & = (\mathrm{holim} \, _N \, B) \wedge_A (\mathrm{holim} \, _N \, B) \overset{?}{\simeq} \mathrm{holim} \, _N \, \mathrm{holim} \, _\mathbf{N} \, (B \wedge_A \mathbf{B}) \\ & \simeq \mathrm{holim} \, _N \, \mathrm{holim} \, _\mathbf{N} \, (\textstyle{\prod}_G \, \mathbf{B}) \cong \mathrm{holim} \, _N \, \textstyle{\prod}_G \, (\mathrm{holim} \, _\mathbf{N} \, \mathbf{B}) \\ & = \mathrm{holim} \, _N \, \textstyle{\prod}_G \, B^{hN} = (\textstyle{\prod}_G \, B^{hN})^{hN} = \bigl(\mathrm{Map}_c(G, B^{hN})\bigr)^{hN}.
\end{align*}
\par
For convenience, we switch to working in $Sp$, where we use homotopy fixed points for a profinite group. Since $\mathrm{Map}_c(G,B^{hN}) \cong \mathrm{Map}_c(N, \prod \, _{G/N} \,B^{hN})$ and $\prod \, _{G/N} \, B^{hN}$ is a fibrant spectrum, $\mathrm{Map}_c(G,B^{hN})$ is fibrant in $Sp_{_N}$. Thus, \begin{align*}\bigl(\mathrm{Map}_c(G, B^{hN})\bigr)^{hN}& = \bigl(\bigl(\mathrm{Map}_c(G, B^{hN})\bigr)_{f,N}\bigr)^{N}\\ & =\bigl(\mathrm{Map}_c(G, B^{hN})\bigr)^{N} = \textstyle{\prod} \, _{G/N} \, B^{hN}.\end{align*}
\end{proof}
\par
We are also interested in an $E$-local version of the notion of Galois extension.
\begin{Def}\label{gallocal}
Let $G$ be a finite group. A map $A \rightarrow B$ of $E$-local commutative $S$-algebras is an {\em $E$-local $G$-Galois extension} if the following conditions hold:
\begin{enumerate}
\item
$G$ acts on $B$ through $A$-algebra maps.
\item
The natural map $A \rightarrow B^{hG}$ is a weak equivalence. 
\item 
There is a weak equivalence $L_E(B \wedge_A B) \cong \prod_G \, B.$ 
\end{enumerate}
\end{Def}
\begin{Rk}
Definition \ref{gallocal} is slightly different from the definition for the same term given in \cite{Rognes2003}, where Rognes does not assume that $A$ and $B$ are $E$-local, but he requires $A \rightarrow B^{hG}$ and $B \wedge_A B \simeq \prod _G B$ to be $E$-equivalences.
\end{Rk}
\par
The result below concerns the map $\alpha(U) = \mathbf{F}(G_n/U \rightarrow G_n/{G_n}) \: E_n^{hG_n} \rightarrow E_n^{hU}.$ Following \cite[Thm. 1(iii)]{DH}, we identify $E_n^{hG_n}$ with $\lhat S^0$.
\begin{Thm}\label{ex1}
The map $\alpha(U) \: \lhat S^0 \rightarrow E_n^{hU}$ of commutative $S$-algebras is a $K(n)$-local $G_n/{U}$-Galois extension.
\end{Thm}
\begin{proof}
As stated in \cite[pg. 8]{LHS}, the finite group $G_n/U{{}}$ acts on $E_n^{hU{{}}}$ through maps of $\lhat S^0$-algebra maps. Also, by \cite[Thm. 4]{DH}, $\lhat S^0 \simeq \bigl(E_n^{hU{{}}}\bigr)^{hG_n/U{{}}}.$ Also, by \cite[Cor. 3.9]{LHS}, there is a weak equivalence $$\lhat(E_n^{hU{{}}} \wedge_{L S^0} E_n^{hU{{}}}) \simeq \mathrm{Map}_c(G_n/{U{{}}}, E_n^{hU{{}}}) \simeq \textstyle{\prod} \, _{G_n/{U{{}}}} \, E_n^{hU{{}}}.$$
\end{proof}
\begin{Rk}
Note that (the author believes) $E_n^{hU{{}}} \wedge_{\lhat S^0} E_n^{hU{{}}}$ is not $K(n)$-local (see \cite[VIII, Cor. 3.5]{EKMM}), so there is no weak equivalence $E_n^{hU{{}}} \wedge_{L S^0} E_n^{hU{{}}} \simeq \textstyle{\prod} \, _{G_n/{U{{}}}} \, E_n^{hU{{}}}.$ Therefore, $\alpha(U{{}})$ is not a $G_n/{U{{}}}$-Galois extension.
\end{Rk}
\begin{Rk}
Following \cite[Chp. VIII]{EKMM} and \cite[\S 1]{LHS}, let $L_E^A$ denote Bousfield localization with respect to $E$ for $A$-modules, where $E$ is an $A$-module. Note that the $K(n)$-local spectrum $K(n)$ is a module over $\lhat S^0$, the unit in the $K(n)$-local category. One can define a {\em strongly $E$-local $G$-Galois extension} to be as in Definition \ref{gallocal}, except that in condition (3), the functor $L_E$ is replaced with $L_E^A$. Then $\alpha(U)$ is a strongly $K(n)$-local $G_n/U$-Galois extension. (To verify this, we only need to show that \begin{equation}\label{verify}\zig L_{K(n)}^{\lhat S^0}(E_n^{hU} \wedge_{\lhat S^0} E_n^{hU}) \simeq \lhat(E_n^{hU} \wedge_{\lhat S^0} E_n^{hU}).\end{equation} Let $X$ be an arbitrary $\lhat S^0$-module. Then by \cite[pg. 4]{LHS}, $\lhat X \simeq L_{\lhat S^0 \wedge K(n)}^{\lhat S^0} (X).$ Since $\lhat S^0 \wedge K(n)$ is $K(n)$-local, $$\lhat S^0 \wedge K(n) \simeq \lhat(\lhat S^0 \wedge K(n)) \simeq \lhat(S^0 \wedge K(n)) \simeq K(n).$$ Thus, $\lhat X \simeq L_{K(n)}^{\lhat S^0} (X),$ and (\ref{verify}) is true.)
\end{Rk}
\par
The lemma below implies that $\lhat(E_n^{hU} \wedge_{\lhat S^0} E_n^{hU})$, which is associated to the extension $\alpha(U)$, and $\lhat(E_n^{hU} \wedge E_n^{hU})$ are the same as $S$-modules.
\begin{Lem}\label{lem}
There is a weak equivalence $\lhat(E_n^{hU} \wedge E_n^{hU}) \simeq \textstyle{\prod} \, _{G_n/U} \, E_n^{hU}.$
\end{Lem}
\begin{proof}
The finite product of $K(n)$-local spectra is $K(n)$-local, so both spectra under consideration are $K(n)$-local. Then it suffices to show that there is a weak equivalence $\lhat(E_n^{hU} \wedge E_n^{hU} \wedge E_n) \simeq \lhat((\textstyle{\prod} \, _{G_n/U} \, E_n^{hU}) \wedge E_n),$ and this follows from $\lhat(E_n^{hU} \wedge E_n) \simeq \textstyle{\prod} \, _{G_n/U} \, E_n$ \cite[Cor. 5.5]{DH}.
\end{proof}
\end{section}
\begin{section}{Associative Galois extensions}
\par
In this section, we consider Galois extensions of $S$-algebras that are not necessarily commutative. We use the fact that if $R$ is just an $S$-algebra and $M$ and $N$ are right and left $R$-modules, respectively, then the tensor product $M \wedge_R N$ is still defined, though it need not be an $R$-module \cite[III, Def. 3.1]{EKMM}. 
\par
Recall that $F_n=\mathrm{colim} \, _i \, E_n^{hU_i}$ is a discrete $G_n$-spectrum of simplicial sets. Also, \cite[Def. 1.5, Thm. 3(i)]{DH} shows that $E_n \simeq \lhat (\mathrm{hocolim} \, ^\mathcal{E} _i E_n^{hU_i})$, where $\mathrm{hocolim} \, ^\mathcal{E}$ is the homotopy colimit in the model category $\mathcal{E}$ of commutative $S$-algebras. Furthermore, by \cite[Lem. 6.2]{DH}, $\mathrm{hocolim} \, ^\mathcal{E} _i E_n^{hU_i} \simeq \mathrm{colim} \, _i \, E_n^{hU_i}$, where the colimit is in the category of $S$-modules. Thus, we can regard $F_n$ as a commutative $S$-algebra. 
\par
\begin{Rk}\label{unproven}
In the next two paragraphs, all statements are unproven, except for the statements that are in {\it italics}, which are known to be true. We include the unproven assertions because, if true, they form an integral part of the story of how Galois extensions appear in Lubin-Tate theory, as the rest of this paper shows. Also, the author believes the assertions are probably true, and he has worked on showing that $E_1 \wedge M(p^i)$ can be realized in $Sp^{\, a}_{_{G_1}}$ (see below). We assume the unproven statements are true for the remainder of the paper.
\end{Rk}
\par
For $G$, a profinite group, there is a model category $Sp^{\, c}_{_G}$ of discrete $G$-commutative symmetric ring spectra, that is, $E_\infty$-objects in the category of symmetric spectra of simplicial sets that are also discrete $G$-spectra, such that the $G$-action is by $E_\infty$-maps. Let $Sp^{\, \Sigma}_{_G}$ be the model category of discrete $G$-symmetric spectra, and let $Sp^{\, c}$ be the model category of commutative symmetric ring spectra. Then the forgetful functor $Sp^{\, c}_{_G} \rightarrow Sp^{\, \Sigma}_{_G}$ and the $G$-fixed points functor $(-)^G \: Sp^{\, c}_{_G} \rightarrow Sp^{\, c}$ preserve all weak equivalences and fibrations. Also, if $X \in Sp^{\, c}_{_G}$, then $X^{hG}$ is a commutative symmetric ring spectrum. 
\par
{\it Now we consider what kind of Galois extension arises for $B=F_n \wedge M_I \simeq E_n \wedge M_I.$ It is widely believed that $F_n \wedge M_I$ cannot be a commutative $S$-algebra. However, it is thought that $F_n  \wedge M_I$ is an $S$-algebra, since Andrew Baker proved that the closely related spectra $E(n)/{I_n^k}$ are $S$-algebras \cite{Baker}.} Further, we suppose that $F_n \wedge M_I$ is a discrete $G_n$-symmetric ring spectrum, that is, $F_n \wedge M_I$ is an object in $Sp^{\, a}_{_{G_n}}$, the model category of $A_\infty$-objects in the category of symmetric spectra that are discrete $G_n$-spectra with an action by $A_\infty$-maps. As above, the forgetful functor $Sp^{\, a}_{_G} \rightarrow Sp^{\, \Sigma}_{_G}$ and the $G$-fixed points functor $(-)^G \: Sp^{\, a}_{_G} \rightarrow Sp^{\, a}$ preserve all weak equivalences and fibrations. ($Sp^{\, a}$ is the model category of symmetric ring spectra.) Thus, if $X \in Sp^{\, a}_{_G}$, then $X^{hG}$ is a symmetric ring spectrum. 
\par
The next result is useful for verifying that certain maps are Galois extensions.
\begin{Lem}\label{jardine}
Let $X$ be a discrete $G$-spectrum and let $N$ be an open normal subgroup of $G$. Then there is a weak equivalence $X^{hG} \rightarrow (X^{hN})^{hG/N}$ is a weak equivalence.
\end{Lem}
\begin{proof}
The sheaf of spectra $\mathrm{Hom}_G(-,X_{f,G})$ is a globally fibrant presheaf of spectra. Then \cite[Prop. 6.39]{Jardine} implies that there is a weak equivalence $$X^{hG} \cong \mathrm{Hom}_G(\ast,X_{f,G}) \rightarrow \mathrm{holim} \, _{G/N} \, \mathrm{Hom}_G(G/N, X_{f,G}) \cong \mathrm{holim} \, _{G/N} \, X^{hN},$$ since $X_{f,G}$ is fibrant in $Sp_{_N}$. Note that $(X_{f,G})^N$ is a $G/N$-spectrum. Since $G/N$ is finite and the $G/N$-spectrum $X^{hN}$ is fibrant in $Sp$, $\mathrm{holim} \, _{G/N} \, X^{hN} = (X^{hN})^{hG/N}$.
\end{proof}
\par
Assuming the hypothetical picture discussed above, the map $$\gamma(U,I) \: \lhat M_I \cong (F_n \wedge M_I)^{hG_n} \rightarrow (F_n \wedge M_I)^{hU} \cong E_n^{hU} \wedge M_I$$ is a map of $S$-algebras. By Lemma \ref{jardine}, $\bigl((F_n \wedge M_I)^{hU}\bigr)^{hG_n/U} \simeq (F_n \wedge M_I)^{hG_n}.$ 
\begin{Lem}
There is a weak equivalence of $S$-modules $$(F_n \wedge M_I)^{hU} \wedge_{\lhat M_I} (F_n \wedge M_I)^{hU} \simeq \textstyle{\prod} \, _{G_n/U} \, (F_n \wedge M_I)^{hU}.$$ 
\end{Lem}
\begin{proof}
By Lemma \ref{lem}, there are weak equivalences $$\textstyle{\prod} \, _{G_n/U} \, (F_n \wedge M_I)^{hU} \simeq  \prod \, _{G_n/U} \, (E_n^{hU} \wedge M_I) \simeq E_n^{hU} \wedge E_n^{hU} \wedge M_I.$$ Now we consider the left hand side of the desired weak equivalence. Note that $$(F_n \wedge M_I)^{hU} \cong F_n^{hU} \wedge M_I \cong F_n^{hU} \wedge L_n(M_I) \cong F_n^{hU} \wedge \lhat M_I.$$ Thus, by \cite[III, Prop. 3.6]{EKMM}, we have: \begin{align*} (F_n \wedge M_I)^{hU} \wedge_{\lhat M_I} (F_n \wedge M_I)^{hU} & \cong \bigl(F_n^{hU} \wedge \lhat M_I \bigr) \wedge_{\lhat M_I} (F_n^{hU} \wedge M_I) \\ & \simeq F_n^{hU} \wedge F_n^{hU} \wedge M_I \cong E_n^{hU} \wedge E_n^{hU} \wedge M_I.\end{align*}
\end{proof}
\par
This lemma motivates us to make the following definition. The theorem below follows immediately from the lemma and the definition.
\begin{Def}\label{non}
Let $G$ be a finite group. A map $A \rightarrow B$ of $S$-algebras is an {\em associative $G$-Galois extension} if the following conditions hold:
\begin{enumerate}
\item
$G$ acts on $B$.
\item
There is an isomorphism $A \cong B^{hG}$ in the stable homotopy category. 
\item
There is a weak equivalence $B \wedge_A B \simeq \prod_G \, B,$ where $B$ is a left and a right $A$-module.
\end{enumerate}
\end{Def}
\begin{Thm}\label{assoc}
The map $\gamma(U,I) \: \lhat M_I \rightarrow (F_n \wedge M_I)^{hU} \cong E_n^{hU} \wedge M_I$ of $S$-algebras is an associative $G_n/U$-Galois extension.
\end{Thm}
\begin{Rk}
Without spelling out another definition, it is easy to see that $\gamma(U,I)$ is also a $K(n)$-local associative $G_n/U$-Galois extension, since $\prod \, _{G_n/U} \, (E_n^{hU} \wedge M_I)$ is $K(n)$-local.
\end{Rk}
\end{section}
\begin{section}{Galois extensions for profinite groups}
\par
In this section we consider the notion of $G$-Galois extension for a profinite group $G$. 
\begin{Def}\label{galpro}
Let $G$ be a profinite group, and let $A \rightarrow B$ be a map of commutative $S$-algebras. Also, let $B$ be a discrete $G$-commutative symmetric ring spectrum, that is, $B \in Sp^{\, c}_{_G}$. Then $A \rightarrow B$ is a {\em $G$-Galois extension} if the following conditions hold:
\begin{enumerate}
\item
There is a compatible $G$-action on $B$ that is by $A$-algebra maps.
\item
There is an isomorphism $A \cong B^{hG}$ in the stable homotopy category.
\item
There is a weak equivalence $B \wedge_A B \simeq \mathrm{colim} \, _N \, \prod \, _{G/N} \, B.$
\end{enumerate}
\end{Def}
\begin{Rk}
The $A$-algebra action on the commutative $S$-algebra $B$ is in the world of $S$-modules, whereas the discrete $G$-action is in the world of symmetric spectra of simplicial sets. Since these categories are different, we only ask for compatibility in condition (1) above, instead of requiring that the discrete $G$-action on $B$ be by $A$-algebra maps. For an example of what ``compatible'' means, see the related example mentioned in Remark \ref{list} (1).
\end{Rk}
\begin{Rk}
If $B$ is a spectrum of simplicial sets, then there is an isomorphism $\mathrm{colim} \, _N \, \prod \, _{G/N} \, B \cong \mathrm{Map}_c(G,B)$. We use the former construction in the above definition, since the latter construction, in general, does not give the right spectrum, if $B$ is a spectrum of topological spaces. (If $B$ is an $S$-module and $V$ is a finite dimensional subspace of $\mathbb{R}^\infty$, then $BV$, in general, is not a discrete space, and $\mathrm{Map}_c(G,BV) \neq \mathrm{Map}_c(G,BV_\mathrm{dis})$, where $BV_\mathrm{dis}$ is the set $BV$ with the discrete topology.)
\end{Rk}
\begin{Rk}
Let $G$ be profinite and let $A \rightarrow B$ be a $G$-Galois extension. Now suppose that $G$ is finite. Then it is known that $B^{hG} = (B_{f,G})^G$ and $\mathrm{Map}_G(EG_+,B)$ are weakly equivalent. Also, since $B$ is a discrete $G$-spectrum, $$\mathrm{colim} \, _N \, \textstyle{\prod} \, _{G/N} \, B \cong \mathrm{Map}_c(G,B)=\textstyle{\prod}_G \, B.$$ Thus, the Galois extension satisfies the conditions of Definition \ref{gal}, so that Definition \ref{galpro} includes Definition \ref{gal} as a special case, as desired.
\end{Rk}
\par
We have the following definition for when $A$ and $B$ are only $S$-algebras.
\begin{Def}\label{galpro2}
Let $G$ be a profinite group, and let $A \rightarrow B$ be a map of $S$-algebras. Also, let $B$ be a discrete $G$-symmetric ring spectrum, that is, $B \in Sp^{\, a}_{_G}$. Then $A \rightarrow B$ is an {\em associative $G$-Galois extension} if the following conditions hold:
\begin{enumerate}
\item
There is an isomorphism $A \cong B^{hG}$ in the stable homotopy category.
\item
There is a weak equivalence $B \wedge_A B \simeq \mathrm{colim} \, _N \, \prod \, _{G/N} \, B.$
\end{enumerate}
\end{Def}
\end{section}
\begin{section}{Filtered Galois extensions}
\par
In this section, we introduce the notion of filtered Galois extension, which is essentially what Rognes calls a pro-Galois extension \cite[pg. 1]{Rognes2000}. We reserve use of the prefix ``pro'' for later, when we consider Galois extensions that are inverse limits of Galois extensions.
\begin{Def}\label{direct}
Let $\{A \rightarrow B_\alpha\}_\alpha$ be a direct system of $G_\alpha$-Galois extensions, with $\{G_\alpha\}_\alpha$ an inverse system of finite groups, such that each map $B_\alpha \rightarrow B_{\alpha'}$ is $G_{\alpha'}$-equivariant. Let $G = \mathrm{lim} \, _\alpha \, G_\alpha$ and let $B = \mathrm{colim} \, _\alpha \, B_\alpha$, so that $G$ is a profinite group and $B \in Sp_{_G}.$ Henceforth, whenever we say {\em direct system of $G_\alpha$-Galois extensions}, we are referring to a system with these properties. A {\em direct system of associative $G_\alpha$-Galois extensions} is a direct system of $G_\alpha$-Galois extensions, except we only require the Galois extensions to be associative. 
\end{Def}
\par
Let $\{A \rightarrow B_\alpha\}_\alpha$ be a direct system of $G_\alpha$-Galois extensions, such that $G$ has finite virtual cohomological dimension. Recall that if $K$ is profinite with $\mathrm{vcd}(K) < \infty$, then, if $Z \in Sp_{_K}$, $Z^{hK} \simeq \mathrm{holim} \, _\Delta \, \bigl(\Gamma^\bullet_K (Z_{f,K})\bigr)^K.$ Then there are conditionally convergent descent spectral sequences $$E_2^{s,t}(\alpha) = H^s(G_\alpha;\pi_t(B_\alpha)) \Rightarrow \pi_{t-s}(B_\alpha^{hG_\alpha}),$$ and $$E_2^{s,t} = H^s_c(G ;\pi_t(B)) \Rightarrow \pi_{t-s}(B^{hG}).$$ Taking a colimit of the spectral sequences $E_r^{\ast,\ast}(\alpha)$ yields the spectral sequence $$\mathrm{colim} \, _\alpha \, E_2^{s,t}(\alpha) \cong H^s_c(G;\pi_t(B)) \Rightarrow \pi_{t-s}\bigl(\mathrm{holim} \, _\Delta \, \mathrm{colim} \, _\alpha \, \bigl(\Gamma^\bullet_{G_\alpha} ((B_\alpha)_{f,G_\alpha})\bigr)^{G_\alpha}\bigr).$$ Thus, the isomorphism of spectral sequences $\mathrm{colim} \, _\alpha \, E_r^{\ast,\ast}(\alpha) \cong E_r^{\ast,\ast},$ for $r \geq 2,$ implies that \begin{equation}\label{neato}\zig B^{hG} \cong \mathrm{holim} \, _\Delta \, \mathrm{colim} \, _\alpha \, \bigl(\Gamma^\bullet_{G_\alpha} ((B_\alpha)_{f,G_\alpha})\bigr)^{G_\alpha}.\end{equation}
\par
Observe that if, in (\ref{neato}), the colimit and the holim commute with each other (that is, if the spectral sequence $\mathrm{colim} \, _\alpha \, E_r^{\ast,\ast}(\alpha)$ converges to the colimit of the abutments $\pi_\ast(B_\alpha^{hG_\alpha})$), then $$B^{hG} \cong \mathrm{colim} \, _\alpha \, B_\alpha^{hG_\alpha} \simeq \mathrm{colim} \, _\alpha \, A = A.$$ However, a strong hypothesis (e.g. the collection $\{E_2^{\ast, \ast}(\alpha)\}$ is uniformly bounded on the right - see \cite[Lem. 5.50]{Thomason}) is  needed for this to be true. Thus, in general, we believe that it need not be the case that, given a directed system $\{A \rightarrow B_\alpha\}_\alpha$ of $G_\alpha$-Galois extensions, there is a weak equivalence $B^{hG} \simeq A.$ Thus, $A \rightarrow B$ in general, is not automatically a $G$-Galois extension. This motivates the following definition. 
\begin{Def}\label{galfil}
Let $\{A \rightarrow B_\alpha\}_\alpha$ be a  direct system of (associative) $G_\alpha$-Galois extensions. As before,  $G = \mathrm{lim} \, _\alpha \, G_\alpha$ is profinite and $B = \mathrm{colim} \, _\alpha \, B_\alpha \in Sp_{_G}.$ If $A \rightarrow B$ is a (associative) $G$-Galois extension, then $A \rightarrow B$ is called a {\em (associative) filtered $G$-Galois extension}. 
\end{Def}
\par
Recall that in Theorem \ref{assoc}, we showed that $\{\gamma(U_i,I)\}_i = \{\lhat M_I \rightarrow E_n^{hU_i} \wedge M_I\}_i$ is a direct system of associative $G_n/{U_i}$-Galois extensions. 
\begin{Thm}\label{thm}
The map $\gamma(I) = \mathrm{colim} \, _i \, \gamma(U_i,I) \: \lhat M_I \rightarrow F_n \wedge M_I$ is an associative filtered $G_n$-Galois extension.
\end{Thm}
\begin{proof}[Proof of Theorem \ref{thm}]
We only have to show that $\gamma(I)$ is an associative $G_n$-Galois extension. Since $(F_n \wedge M_I)^{hG_n} \cong E_n^{hG_n} \wedge M_I \cong \lhat M_I,$ it suffices to show that $(F_n \wedge M_I) \wedge_{\lhat M_I} (F_n \wedge M_I) \simeq \mathrm{colim} \, _i \, \textstyle{\prod} \, _{G_n/{U_i}} \, (F_n \wedge M_I).$ Since $F_n$ is $E(n)$-local, $F_n \wedge M_I \simeq F_n \wedge \lhat M_I,$ so that \begin{align*}(F_n \wedge M_I) \wedge_{\lhat M_I} (F_n \wedge M_I) & \simeq (F_n \wedge \lhat M_I) \wedge_{\lhat M_I} (F_n \wedge M_I) \simeq F_n \wedge F_n \wedge M_I \\ & \simeq E_n \wedge E_n \wedge M_I \simeq \mathrm{Map}_c(G_n, F_n \wedge M_I) \\ & \simeq \mathrm{colim} \, _i \, \textstyle{\prod} \, _{G_n/{U_i}} \, (F_n \wedge M_I).\end{align*}
\end{proof}
\begin{Rk}
The map $\gamma(I)$ is also a $K(n)$-local associative $G_n$-Galois extension, where we use the following definition.
\end{Rk}
\begin{Def}
Let $G$ be a profinite group, and let $A \rightarrow B$ be a map of $E$-local $S$-algebras. Also, let $B$ be a discrete $G$-symmetric ring spectrum, that is, $B \in Sp^{\, a}_{_G}$. Then $A \rightarrow B$ is an {\em $E$-local associative $G$-Galois extension} if the following conditions hold:
\begin{enumerate}
\item
There is an isomorphism $A \cong B^{hG}$ in the stable homotopy category.
\item
There is a weak equivalence $L_E(B \wedge_A B) \simeq \mathrm{colim} \, _N \, \prod \, _{G/N} \, B.$
\end{enumerate}
\end{Def}
\par
Let (a) $\{A \rightarrow B_\alpha\}_\alpha$ be a direct system of (associative) $G_\alpha$-Galois extensions, and (b) assume that $A \rightarrow B$ is a map of commutative $S$-algebras. Note that if $X$ is a discrete $G$-set such that $X \cong \mathrm{colim} \, _N \, X(N)$, where each $X(N)$ is a $G/N$-set, then $\mathrm{Map}_c(G,X) \cong \mathrm{colim} \, _N \, \mathrm{Map}_c(G/N, X(N))$ (see e.g. \cite[Lem. 6.5.4(a)]{Ribes}). Similarly, \begin{align*} \mathrm{Map}_c(G,B) & \cong \mathrm{Map}_c(\mathrm{lim} \, _\alpha\, G_\alpha, \mathrm{colim} \, _\alpha \, B_\alpha) \cong \mathrm{colim} \, _\alpha \, \mathrm{Map}_c(G_\alpha, B_\alpha) \\ & \cong \mathrm{colim} \, _\alpha \, \textstyle{\prod} \, _{G_\alpha} \, B_\alpha \simeq \mathrm{colim} \, _\alpha \, (B_\alpha \wedge_A B_\alpha).
\end{align*}
\par
Let $\{\alpha'\}$ be a copy of the indexing set $\{\alpha\}$, so that $\alpha = \alpha'$. Then the set of pairs $\{(\alpha, \alpha)\}_\alpha$ is cofinal in the indexing set $\{(\alpha, \alpha')\}_{\alpha, \alpha'}$ of all pairs, so that $$\mathrm{colim} \, _\alpha \, (B_\alpha \wedge_A B_\alpha) \cong \mathrm{colim} \, _{(\alpha, \alpha)} \, (B_\alpha \wedge_A B_\alpha) \cong \mathrm{colim} \, _{(\alpha, \alpha')} \, (B_\alpha \wedge_A B_{\alpha'}).$$ Since the construction $B_\alpha \wedge_A B_\alpha$ is a coequalizer, \begin{align*}\mathrm{colim} \, _{(\alpha, \alpha')} \, (B_\alpha \wedge_A B_{\alpha'}) & \cong (\mathrm{colim} \, _\alpha \, B_\alpha) \wedge_A (\mathrm{colim} \, _{\alpha'} \, B_{\alpha'}) \\ & \cong B \wedge_A B.\end{align*} Thus, $\mathrm{Map}_c(G,B) \simeq B \wedge_A B$, and we summarize this discussion in the remark below.
\begin{Rk}
As stated in \cite[pg. 1]{Rognes2000}, (a) and (b) above are enough to imply the weak equivalence $B \wedge_A B \simeq \mathrm{colim} \, _N \, \prod \, _{G/N} \, B.$ Thus, Definition \ref{galfil} can be simplified by noting that the last condition in Definitions \ref{galpro} and \ref{galpro2} can be ignored.
\end{Rk}
\end{section}
\begin{section}{A consequence of Theorem \ref{unsure}, when $G$ is profinite}
In this brief section, we assume that Theorem \ref{unsure} is true, when $G$ is profinite. Thus, we are assuming that if (i) $G$ is profinite; (ii) $A \rightarrow B$ is a $G$-Galois extension of commutative $S$-algebras; and (iii) $N$ is an open normal subgroup of $G$, then $A \rightarrow B^{hN}$ is a $G/N$-Galois extension.
\begin{Rk}\label{remark} 
In the theorem below, $N$ is an open normal subgroup of $G$, $B$ is a discrete $G$-commutative symmetric ring spectrum, and $B_{f,G}$ comes from factoring $B \rightarrow \ast$ in $Sp^{\, c}_{_G}$, as $B \rightarrow B_{f,G} \rightarrow \ast$, a trivial cofibration followed by a fibration. Since the forgetful functor $Sp^{\, c}_{_G} \rightarrow Sp^{\, \Sigma}_{_G}$ preserves weak equivalences and fibrations, $B \rightarrow B_{f,G}$ is a weak equivalence in $Sp^{\, \Sigma}_{_G}$, and $B_{f,G}$ is fibrant in $Sp^{\, \Sigma}_{_G}$. Thus, in $Sp^\Sigma_{_N}$, $B \rightarrow B_{f,G}$ is a weak equivalence and $B_{f,G}$ is fibrant, so that $(B_{f,G})^N$ is a model for $B^{hN}$.
\end{Rk}
\begin{Thm}\label{canonical}
Let $G \cong \mathrm{lim} \, _N \, G/N$ be profinite. If $\, A \rightarrow B$ is a $G$-Galois extension, then the direct system $\{A \rightarrow B^{hN}\}_N$, of $G/N$-Galois extensions, makes the map $A \rightarrow B_{f,G}$ a filtered $G$-Galois extension in a canonical way.
\end{Thm}
\begin{proof}
Observe that $B_{f,G} = \bigcup \, _N \, (B_{f,G})^N = \mathrm{colim} \, _N \, B^{hN},$ as required. Also, $(B_{f,G})^{hG} \simeq B^{hG} \simeq A,$ and $(B_{f,G}) \wedge_A (B_{f,G}) \simeq B \wedge_A B \simeq \mathrm{colim} \, _N \, \prod \, _{G/N} \, B \simeq \mathrm{colim} \, _N \, \prod \, _{G/N} \, B_{f,G}.$
\end{proof}
\begin{Rk}
This theorem says that every $G$-Galois extension is canonically a filtered $G$-Galois extension.
\end{Rk}
\end{section} 
\begin{section}{Pro-Galois extensions}
In this section, we define a notion of Galois extension for towers of discrete $G$-spectra. We are primarily interested in understanding the structure of the Galois extension $\lhat S^0 \rightarrow E_n$, which Rognes has referred to as a ``$K(n)$-local $G_n$-pro-Galois extension'' \cite{Rognes2003}. We begin by recalling that $\gamma(I) \: \lhat M_I \rightarrow F_n \wedge M_I \overset{\simeq}{\longrightarrow} E_n \wedge M_I$ is an associative filtered $G_n$-Galois extension. Thus, $\{\gamma(I)\}_I$ is an inverse system of associative filtered $G_n$-Galois extensions.
\begin{Rk}
The definition below is only in the $K(n)$-local setting because this is all that is needed for our examples.
\end{Rk}
\begin{Def}\label{pro}
Let $J= \{\cdots \rightarrow i \rightarrow i-1 \rightarrow \cdots \rightarrow 1 \rightarrow 0\}.$
Let $\{A_i \rightarrow B_i\}_i$ be a $J$-shaped tower of $K(n)$-local $G$-Galois extensions, such that $\{B_i\}$ is a tower in $Sp_{_G}$, and the isomorphism $B_i^{hG} \cong A_i$ comes from a natural weak equivalence $A_i \rightarrow B_i^{hG}$. (Whenever $B_i$ is viewed as an object of $Sp_{_G}$, then it is assumed to be fibrant there. Similarly, whenever $A_i$ is viewed as an object of $Sp$, then it is assumed to be fibrant.) We allow any or all of the extensions $A_i \rightarrow B_i$ to be associative. Let $A = \mathrm{holim} \, _i \, A_i$, $B = \mathrm{holim} \, _i \, B_i$, and let $A \rightarrow B$ be the obvious map. Then $A \rightarrow B$ is a {\em $K(n)$-local $G$-pro-Galois extension} if the following conditions hold:
\begin{enumerate}
\item
The map $A \rightarrow B$ is a map of commutative $S$-algebras.
\item
The spectrum $B$ is a continuous $G$-spectrum, $G$ acts on $B$ by maps of $A$-algebras, and these two $G$-actions are compatible.
\item
There is a weak equivalence $A \simeq B^{hG}.$
\item
There is a weak equivalence $\lhat(B \wedge_A B) \simeq \mathrm{holim} \, _i \, \bigl(\mathrm{colim} \, _N \, \prod \, _{G/N} \, B_i\bigr).$
\end{enumerate}
\end{Def}              
\begin{Rk}\label{list} 
We make some remarks about this definition; in particular, we discuss what is and is not automatically entailed by the hypotheses of the definition.
\begin{enumerate}
\item
So that condition (2) above is actually met in practice, we do not require that the continuous action be by maps of $A$-algebras; we only require that the continuous action and the $A$-algebra action be compatible. For example, $G_n$ acts on $E_n$ by maps of $\lhat S^0$-algebras and this action yields the continuous action described in \cite{thesis}, but the continuous action is only (thus far, known to be) by maps of (unstructured) spectra.
\item
Since all the $A_i$ and $B_i$ are $K(n)$-local,
the homotopy limits $A$ and $B$ are also $K(n)$-local.
\item
The hypotheses of the definition imply that $B$ is automatically a continuous $G$-spectrum.
\item
By \cite{thesis}, $B^{hG} = \mathrm{holim} \, _i \, B_i^{hG} \overset{\simeq}{\longleftarrow} \mathrm{holim} \, _i \, A_i = A,$ so that the assumptions automatically imply that condition (3) holds. 
\end{enumerate}
\end{Rk}
\begin{Rk}
We explain part of our motivation for condition (4) in Definition \ref{pro}. Recall (from \cite{thesis}) that the functor $\mathrm{Map}_c(G, -) \: Sp \rightarrow Sp_{_G}$ is a right Quillen functor. Let $X \in Sp_{_G}$ be fibrant, so that $X$ is also fibrant in $Sp$, and hence, $\mathrm{colim} \, _N \, \prod \, _{G/N} \, X \cong \mathrm{Map}_c(G,X)$ is fibrant in $Sp_{_G}$. Then $\, \cdots \rightarrow X \rightarrow X$, the constant tower of fibrations of fibrant spectra in $Sp$, gives \begin{align*}\mathrm{holim} \, _i \, \bigl(\mathrm{colim} \, _N \, \textstyle{\prod} \, _{G/N} \, X \bigr) & \cong \mathrm{holim} \, _i \, \mathrm{Map}_c(G,X) \simeq \mathrm{lim} \, _i \, \mathrm{Map}_c(G,X) \\ & \cong \mathrm{Map}_c(G,X) \cong \mathrm{colim} \, _N \, \textstyle{\prod} \, _{G/N} \, X,\end{align*} where the last spectrum, as desired, has the form of the right-hand side in condition (3) of Definition \ref{galpro}. Therefore, a $K(n)$-local $G$-pro-Galois extension is a generalization of a $K(n)$-local (associative) $G$-Galois extension from the setting of discrete $G$-spectra to that of towers of discrete $G$-spectra. 
\end{Rk}
\begin{Thm}\label{ex4}
The map of commutative $S$-algebras $\mathrm{holim} \, _I \, \gamma(I),$ $$\lhat S^0 \cong \mathrm{holim} \, _I \, \lhat M_I \rightarrow \mathrm{holim} \, _I \, (E_n \wedge M_I) \cong E_n,$$ is a $K(n)$-local $G_n$-pro-Galois extension.
\end{Thm}
\begin{proof}
We only need to verify condition (4) of Definition \ref{pro}: by \cite[Cor. 3.9]{LHS}, $$\lhat(E_n \wedge_{\lhat S^0} E_n) \simeq \mathrm{holim} \, _I \, \mathrm{Map}_c(G_n, F_n \wedge M_I) \simeq \mathrm{holim} \, _I \, \mathrm{colim} \, _i \, \textstyle{\prod} \, _{G_n/{U_i}} \, (E_n \wedge M_I).$$ 
\end{proof}
\begin{Rk}
Since $L_{K(n)}^{\lhat S^0}(E_n \wedge_{\lhat S^0} E_n) \simeq \lhat(E_n \wedge_{\lhat S^0} E_n)$, $\lhat S^0 \rightarrow E_n$ is a strongly $K(n)$-local $G_n$-pro-Galois extension.
\end{Rk}
\par
Though we have shown that $\lhat S^0 \rightarrow E_n$ is a $K(n)$-local $G_n$-pro-Galois extension, this notion still does not capture all of the structure that is present in this extension, due to the extra structure that comes from the filtered extension $\gamma(I)$. We capture this additional structure in the following way.
\par
Let $\{A \overset{f_\alpha}{\longrightarrow} B_\alpha\}_\alpha$ be a direct system of, possibly $K(n)$-local, $G_\alpha$-Galois extensions, with $\{G_\alpha\}_\alpha$ an inverse system of finite groups, such that each map $B_\alpha \rightarrow B_{\alpha'}$ is $G_{\alpha'}$-equivariant. As usual, let $G = \mathrm{lim} \, _\alpha \, G_\alpha$, and let $B = \lhat(\mathrm{colim} \, _\alpha \, B_\alpha)$. If $B$ is regarded as a spectrum of simplicial sets, then, letting $(-)_\mathtt{f}$ denote functorial fibrant replacement in $Sp$, $B \cong \mathrm{holim} \, _I \, \mathrm{colim} \, _\alpha \, (B_\alpha \wedge L_n M_I)_\mathtt{f}$ is a continuous $G$-spectrum, since $\mathrm{colim} \, _\alpha \, (B_\alpha \wedge L_n M_I)_\mathtt{f}$ is a discrete $G$-spectrum that is fibrant in $Sp$. Then, if $\lhat(\mathrm{colim} \, _\alpha \, f_\alpha) \: \lhat A \rightarrow B$ is a $G$-Galois extension, we call $\lhat(\mathrm{colim} \, _\alpha \, f_\alpha)$ a {\em $K(n)$-local filtered $G$-Galois extension}. 
\par
Since the direct system $\{\alpha(U_i) \: \lhat S^0 \rightarrow E_n^{hU_i}\}_i$, of $K(n)$-local $G_n/{U_i}$-Galois extensions, yields the extension $\lhat(\mathrm{colim} \, _i \, \alpha(U_i)) \: \lhat S^0 \rightarrow E_n$, we can refer to the map $\lhat S^0 \rightarrow E_n$ as a {\em $K(n)$-local filtered $G_n$-pro-Galois extension.}
\end{section}
\begin{section}{More examples, for closed subgroups of $G_n$}
In this section, for any closed subgroup $G$ of $G_n$ (so $G$ is always profinite and not necessarily finite), we give two examples of Galois extensions. First of all, we slightly expand the definition of filtered $G$-Galois extension.
\begin{Def}
Let $\{A_\alpha \rightarrow B_\alpha\}_\alpha$ be a direct system of (associative) $G_\alpha$-Galois extensions, with $\{G_\alpha\}$ an inverse system of finite groups, where each map $B_\alpha \rightarrow B_{\alpha'}$ is $G_{\alpha'}$-equivariant. Let $G$ and $B$ be defined as usual, and let $A = \mathrm{colim} \, _\alpha \, A_\alpha.$ If $A \rightarrow B$ is a (associative) $G$-Galois extension, then $A \rightarrow B$ is a {\em (associative) filtered $G$-Galois extension}. 
\end{Def}
\begin{Thm}\label{ex5}
The direct system $\{\beta(G,i,I) \: E_n^{hU_iG} \wedge M_I \rightarrow E_n^{hU_i} \wedge M_I\}_i,$ of associative $U_iG/{U_i}$-Galois extensions, makes the map $$\beta(G,I) = \mathrm{colim} \, _i \, \beta(G,i,I) \: (E_n^{hG} \wedge M_I) \cong (F_n \wedge M_I)^{hG} \rightarrow F_n \wedge M_I \simeq E_n \wedge M_I$$ an associative filtered $G$-Galois extension.
\end{Thm}
\begin{proof}
To make the notation more manageable, we use $X/I$ to denote the spectrum $X \wedge M_I$. Since (we are assuming that) $F_n /I \in Sp^{\, a}_{_{G_n}}$, $(F_n /I)^{hU_i} \cong E_n^{hU_i} /I$ and, similarly, $E_n^{hU_iG} /I$ are $S$-algebras. Now we show that 
\begin{equation}\label{interesting}\zig 
(E_n^{hU_i} /I) \wedge_{(E_n^{hU_iG} /I)} (E_n^{hU_i} /I) \simeq \textstyle{\prod} \, _{U_iG/{U_i}} \, (E_n^{hU_i} /I).\end{equation}
\par
Note that $E_n^{hU_i} /I$ and $\textstyle{\prod} \, _{U_iG/{U_i}} \, (E_n^{hU_iG} /I)$ are $K(n)$-local. Applying \cite[Cor. 5.5]{DH}, $$\pi_\ast(E_n \wedge E_n^{hU_i} /I) \cong \pi_\ast(\textstyle {\prod} \, _{G_n/{U_i}} \, (E_n /I)) \cong \textstyle {\prod} \, _{G_n/{U_i}} \, \pi_\ast(E_n /I).$$ Similarly, \begin{align*}\pi_\ast(E_n \wedge (\textstyle{\prod} \, _{U_iG/{U_i}} \, (E_n^{hU_iG} /I))) & \cong \textstyle{\prod} \, _{G_n/{U_iG}} \, \textstyle{\prod} \, _{U_iG/{U_i}} \, \pi_\ast(E_n /I) \\ & \cong \textstyle{\prod} \, _{G_n/{U_i}} \, \pi_\ast(E_n /I).\end{align*} Thus, $\pi_\ast(E_n \wedge E_n^{hU_i} /I) \cong \pi_\ast(E_n \wedge (\textstyle{\prod} \, _{U_iG/{U_i}} \, (E_n^{hU_iG} /I))),$ showing that $$E_n^{hU_i} /I \cong \textstyle{\prod} \, _{U_iG/{U_i}} \, (E_n^{hU_iG} /I).$$ 
\par
This implies that \begin{align*}(E_n^{hU_i} /I) \wedge_{(E_n^{hU_iG} /I)} (E_n^{hU_i} /I) & \cong \bigl(\textstyle{\prod} \, _{U_iG/{U_i}} \, (E_n^{hU_iG} /I)\bigr)\wedge_{(E_n^{hU_iG} /I)} (E_n^{hU_i} /I) \\ & \cong \textstyle{\prod} \, _{U_iG/{U_i}} \, \Bigl((E_n^{hU_iG} /I) \wedge_{(E_n^{hU_iG} /I)} (E_n^{hU_i} /I)\Bigr) \\ & \simeq  \textstyle{\prod} \, _{U_iG/{U_i}} \, (E_n^{hU_i} /I),\end{align*} verifying (\ref{interesting}). This shows that $\beta(G,i,I)$ is an associative $U_iG/{U_i}$-Galois extension.
\par
Note that $E_n/I$ and $\mathrm{colim} \, _i \, \textstyle{\prod} \, _{G/{(U_i \cap G)}} \, (E_n^{hG} /I)$ are $K(n)$-local, and there is an isomorphism $\pi_\ast(E_n \wedge E_n/I) \cong \mathrm{Map}_c(G_n, \pi_\ast(E_n/I)).$ Also, as abelian groups, 
$$\pi_\ast(E_n \wedge (\mathrm{colim} \, _i \, \textstyle{\prod} \, _{G/{(U_i \cap G)}} \, (E_n^{hG} /I))) \cong \mathrm{colim} \, _i \, \textstyle{\prod} \, _{G/{(U_i \cap G)}} \, \pi_\ast(E_n \wedge E_n^{hG}/I),$$ which, by \cite[Prop. 6.3]{DH}, is isomorphic to $$\mathrm{colim} \, _i \, \textstyle{\prod} \, _{G/{(U_i \cap G)}} \, \mathrm{colim} \, _j \, \mathrm{Map}_c(G_n/{U_jG},\pi_\ast(E_n/I)).$$ This last abelian group is isomorphic to 
\begin{align*}
\mathrm{colim} \, _i \, \textstyle{\prod} \, _{G/{(U_i \cap G)}} \, \mathrm{Map}_c(G_n/{G},\pi_\ast(E_n/I)) & \cong \mathrm{Map}_c(G \times G_n/G, \pi_\ast(E_n/I)) \\ & \cong \mathrm{Map}_c(G_n, \pi_\ast(E_n/I)).
\end{align*}
Thus, $E_n/I \cong \mathrm{colim} \, _i \, \textstyle{\prod} \, _{G/{(U_i \cap G)}} \, (E_n^{hG}/I),$ and therefore,
\begin{align*}(E_n/I) \wedge_{(E_n^{hG}/I)} (E_n/I) & \cong \mathrm{colim} \, _i \, \textstyle{\prod} \, _{G/{(U_i \cap G)}} \, \bigl((E_n^{hG}/I) 
\wedge_{({E_n^{hG}}/I)} (E_n/I)\bigr) \\ & \cong \mathrm{colim} \, _i \, \textstyle{\prod} \, _{G/{(U_i \cap G)}} \, (E_n/I),\end{align*} completing the proof.
\end{proof}
\begin{Thm}\label{ex6}
The inverse system $\{\beta(G,I)\}_I$ of associative $G$-Galois extensions makes
the map $$\beta(G) = \mathrm{holim} \, _I \, \beta(G,I) \: E_n^{hG} \rightarrow E_n$$ a $K(n)$-local $G$-pro-Galois extension.
\end{Thm}
\begin{proof}
Using the preceding theorem, it is easy to see that each $\beta(G,I)$ is a $K(n)$-local associative $G$-Galois extension, since $E_n \wedge M_I, E_n^{hG} \wedge M_I,$ and $$\mathrm{colim} \, _i \, \textstyle{\prod} \, _{G/{(U_i \cap G)}} \, (E_n \wedge M_I) \cong (\mathrm{colim} \, _i \, \textstyle{\prod} \, _{G/{(U_i \cap G)}} \, E_n) \wedge M_I$$ are all $K(n)$-local. 
\par   
By \cite[Cor. 3.9]{LHS}, $$\pi_\ast(\lhat(E_n \wedge_{E_n^{hG}} E_n)) \cong \mathrm{Map}_c(G, \pi_\ast(E_n)) \cong \mathrm{lim} \, _I \, \mathrm{Map}_c(G, \pi_\ast(E_n \wedge M_I)).$$ This implies that
\begin{align*}\lhat(E_n \wedge_{E_n^{hG}} E_n) & \cong \mathrm{holim} \, _I \, \mathrm{Map}_c(G, (F_n \wedge M_I)_{f,G}) \\ & \cong \mathrm{holim} \, _I \, \mathrm{colim} \, _i \, \textstyle{\prod} \, _{G/{(U_i \cap G)}} \, (F_n \wedge M_I),\end{align*} where the second expression only occurs in $Sp.$
\end{proof}
Our last result follows from the last line of the proof of Theorem \ref{ex1}.
\begin{Thm}\label{ex7}
The map $\alpha(U) = \mathrm{holim} \, _I \, \gamma(U,I)$ is a $K(n)$-local $G_n/U$-pro-Galois extension.
\end{Thm}
\end{section}

\bibliographystyle{plain}
\bibliography{biblio}
\end{document}